\documentclass[11pt]{article}
\usepackage{amssymb,amsfonts,amsmath,latexsym, epsfig,mathrsfs,colordvi, xcolor, graphicx, cases}
\newtheorem{theo}{Theorem}

\newtheorem{exam}[theo]{Example}
\newtheorem{lem} [theo]{Lemma}
\newtheorem{coro}[theo]{Corollary}

\newtheorem{rem}[theo]{Remark}

\makeatletter \@addtoreset{equation}{section}
\@addtoreset{theo}{}\makeatother

\def\qed{\hfill \rule{4pt}{7pt}}
\def\pf{\noindent {\it Proof.} }

\begin{document}

\title{Counting vertices with given outdegree in plane trees and $k$-ary trees }

\author{Rosena R. X. Du\footnote{Corresponding Author. Email: rxdu@math.ecnu.edu.cn.}, Jia He and Xueli Yun\\ \\ Department of Mathematics, Shanghai Key Laboratory of PMMP \\East China Normal University,
500 Dongchuan Road \\Shanghai, 200241, P. R. China.}

\date{Jan 29, 2018}

\maketitle
\noindent {\bf Abstract:}
We count the number of vertices with given outdegree in plane trees and $k$-ary trees, and get the following results: the total number of vertices of outdegree $i$ among all plane trees with $n$ edges is ${2n-i-1 \choose n-1}$;  the total number of vertices of degree $i$ among all plane trees with $n$ edges is twice this number; and the total number of vertices of outdegree $i$ among all $k$-ary trees with $n$ edges is ${k\choose i}{kn\choose n-i}$. For all these results we bijective proofs.



\section{Introduction}
In this paper we
will assume all the trees are ({\em unlabelled}) {\em plane trees}, i.e., rooted trees whose vertices are considered to be
indistinguishable, but the subtrees at any vertex are linearly
ordered. Let $\mathcal{T}_{n}$ be the set of all \emph{$n$-plane trees}, i.e., plane trees with $n$ edges. For each vertex $v$, we say that $v$ is of {\em
outdegree} $i$ if it has $i$ subtrees, and call vertices of outdegree
$0$ {\em leaves}. Vertices that are not leaves are called {\em internal vertices}.  An ({\em ordinary}) {\em $k$-ary tree} is a plane tree where each vertex has exactly $k$
ordered subtrees, which could be empty. A {\em complete $k$-ary tree} is a $k$-ary tree
for which each internal vertex has exactly $k$ nonempty subtrees. In other words, a complete $k$-ary tree
is a plane tree each of whose internal vertices has outdegree $k$.

%

The set of plane trees is one of the most well-known and well-studied combinatorial structures. Existing results have focused on studying various statistics on plane trees and finding bijections between plane trees and other structures. The enumeration of plane trees by outdegree sequences has been well studied in the literature; see \cite{degreeseq20,degreeseq22,degreeseq54}. 

Counting vertices in plane trees according to the outdegrees was also studied more recently. In \cite{Fine} Deutsch and Shapiro proved that the total number of vertices of odd outdegree over all plane trees with $n$ edges is
\[\frac{2}{3}{2n-1 \choose n}+\frac{1}{3}F_{n-1},\]
where $F_{n}$ is the $n$-th {\em Fine number} (Sequence A000957 in \cite{sequence}).

In this paper we focus on counting the total number of vertices with given outdegree $i$ among all plane trees and $k$-ary trees with $n$-edges. Our first result is the following.

\begin{theo}\label{main_plane}
For integers $n \geq 1$ and $i \geq 0$, the total number of vertices of outdegree $i$ among all plane trees with $n$ edges is \[{2n-i-1 \choose n-1}.\]
\end{theo}

Combining Deutsch and Shapiro’s result with the above theorem we immediately have the following formula for the Fine numbers:

\begin{equation}
F_n=3\sum_{k\geq 0} {2n-2k \choose n}-2{2n+1 \choose n}.
\end{equation}

Deutsch and Shapiro also counted the number of vertices of odd {\em degree} in \cite{Fine}. Here the {\em degree} of a vertex $v$ means the total number of vertices adjacent to $v$. Deutsch and Shapiro proved that among all plane trees with $n$ edges, the total number of vertices of odd degree is twice the total number of vertices of odd outdegree. Their proof is based on a generating function method. Later in \cite{oddeven}, Eu, Liu and Yeh gave a combinatorial proof of this result by defining a two-to-one correspondence. 

Although the authors didn't point out this
 explicitly in \cite{oddeven}, their correspondence showed that over all plane trees with $n$ edges, the number of vertices of degree $i$ is twice the number of vertices of outdegree $i$. Thus combining Eu, Liu and Yeh’s result with Theorem \ref{main_plane} we get the following result.

\begin{coro}\label{coro_plane}
For integers $n \geq 1$ and $i \geq 1$, the total number of vertices of degree $i$ among all plane trees with $n$ edges is \[2{2n-i-1 \choose n-1}.\]
\end{coro}

We make a different approach to this intriguing relation between number of vertices of degree $i$ and number of vertices of outdegree $i$ among all $n$-plane trees.  
Note that for each vertex $v$ in a plane tree with outdegree $i$, the degree of $v$ is also $i$ when $v$ is the root, otherwise the degree of $v$ is $i+1$. We count the number of $n$-plane trees with specified root degree, and the number of vertices with given outdegree among all $n$-plane trees with specified root degree, and therefore prove Corollary \ref{coro_plane} directly.

Another main result of this paper is counting vertices with given outdegree among all $k$-ary trees with $n$ edges. By applying the bijection we define for $n$-plane trees  to $k$-ary trees, we transform this problem to the enumeration of compositions of integers with certain restrictions, and we get the following result. 

\begin{theo}\label{main_kary}
For integers $n,k \geq 1$ and $i \geq 0$, the total number of vertices of outdegree $i$ among all $k$-ary trees with $n$ edges is \[{k\choose i}{kn\choose n-i}.\]
\end{theo}

The rest of the paper is organized as follows. In Section 2 we give a bijective proof for Theorem \ref{main_plane}.  In Section 3, we first give a bijective proof for Theorem \ref{rootdeg} and then prove Corollary \ref{coro_plane}. In Section 4 we give a bijective proof for Theorem \ref{main_kary}. 

\begin{rem}
Gernerating function proofs for Theorem \ref{main_plane} and \ref{main_kary}  can be found in an earlier version of this paper \cite{DuHeYun}. 
\end{rem}

\section{Counting vertices with given outdegree in plane trees}

Let $\mathcal{T}_{n,i}$ denote the set of ordered pairs $(T,v)$ such that $T \in \mathcal{T}_n$, and $v$ is a vertex of $T$ of outdegree $i$. We will associate with each pair $(T,v)\in \mathcal{T}_{n,i}$ a {\em composition} of nonnegative integers $n$. Here a composition of $n$ can be thought of as an expression of $n$ as an {\em ordered} sum of nonnegative integers. More precisely, a sequence $\alpha=(a_1,\ldots, a_k)$ of nonnegative integers satisfying $\sum_{j=1}^{k} a_j=n$ is called a {\em $k$-composition of $n$}. We use $\mathcal{A}_{n,k}$ to denote the set of all $k$-compositions of $n$.

For each $\alpha \in \mathcal{A}_{n,k}$ we define $f(\alpha)=\sum_{j=1}^{k}(a_j -1)$. If a $k$-composition $\alpha=(a_1,\ldots, a_k)$ satisfies $f(a_1, \ldots, a_j)\geq 0$ for each $j, 1\leq j \leq k-1$ and $f(\alpha)=-1$, then we call $\alpha$ a {\em unit composition}. If $\alpha=(a_1,\ldots, a_k)$ satisfies $f(a_1, \ldots, a_j)\geq 0$ for each $j, 1\leq j \leq k$, then we call $\alpha$ a {\em positive composition}. We use $\mathcal{B}$ to denote the set of all unit compositions and $\mathcal{B}_{n,k}$ to denote the set of all unit $k$-compositions of integer $n$.

For each $n$-plane tree $T$, we label the $n+1$ vertices of $T$ with numbers $1$ to $n+1$ in the {\em depth-first order}, or {\em preorder}. (The definition of the depth-first order can be found in \cite[p.336]{Knuth1} or \cite[p.33]{EC2}.) Now we associate with $T$ a composition $\delta(T)$ by $\delta(T)=(d_1, d_2, \ldots, d_{n+1})$, where $d_j$ is the outdegree of the vertex in $T$ with label $j$.


\begin{lem} \label{planeB}
The map $T \mapsto \delta(T)$ is a bijection between $\mathcal{T}_n$ and $\mathcal{B}_{n,n+1}$.
\end{lem}

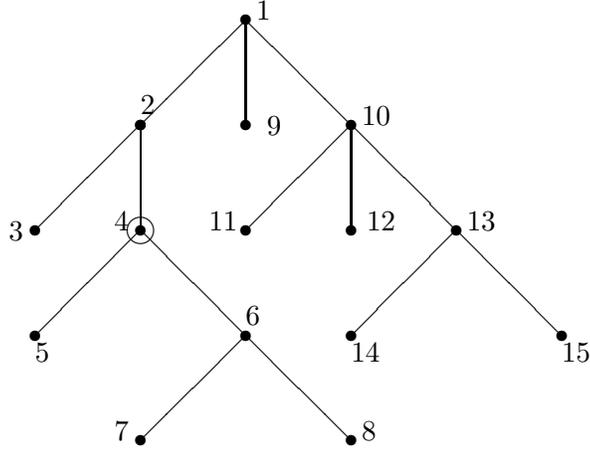
\begin{figure}
\begin{center}
\begin{picture}(220,170)
\setlength{\unitlength}{7mm}

\put(5,8){\line(-1,-1){2}}\put(5,8){\line(1,-1){2}} \put(5,8){\line(0,-1){2}}
\put(5,8){\circle*{0.2}}\put(5.2,8){$1$}

\put(3,6){\circle*{0.2}}\put(3,6.2){$2$}
\put(3,6){\line(-1,-1){2}}\put(1,4){\circle*{0.2}}\put(0.5,3.8){$3$}
\put(3,6){\line(0,-1){2}}\put(3,6){\circle*{0.2}}

\put(5,6){\circle*{0.2}}\put(5.4,5.8){$9$}
\put(7,6){\circle*{0.2}} \put(7.2,6){$10$}
\put(7,6){\line(-1,-1){2}}
\put(7,6){\line(1,-1){2}}
\put(7,6){\line(0,-1){2}}\put(7,4){\circle*{0.2}}\put(7.3,4){$12$}

\put(3,4){\circle*{0.2}}\put(3,4){\circle{0.5}}\put(3,4){\line(-1,-1){2}}
\put(3,4){\line(1,-1){2}}\put(2.5,4){$4$}
\put(5,4){\circle*{0.2}}\put(4.3,4){$11$}
\put(9,4){\circle*{0.2}}\put(9,4){\line(-1,-1){2}}
\put(9,4){\line(1,-1){2}}\put(9.2,4){$13$}

\put(1,2){\circle*{0.2}}\put(1,1.5){$5$}
\put(5,2){\circle*{0.2}}\put(5,2.2){$6$}
\put(5,2){\line(-1,-1){2}}\put(3,0){\circle*{0.2}}\put(2.5,0){$7$}
\put(5,2){\line(1,-1){2}}\put(7,0){\circle*{0.2}}\put(7.2,0){$8$}

\put(7,2){\circle*{0.2}}\put(7,1.5){$14$}
\put(11,2){\circle*{0.2}}\put(11,1.5){$15$}

\end{picture}
\caption{A plane tree with 14 edges, with one of its vertices of outdegree $2$ circled.}\label{Fig_plane}
\end{center}
\end{figure}

\begin{exam}
Let $n=14$. Figure \ref{Fig_plane} shows a plane tree $T$ with 14 edges, whose vertices are labeled in depth-first order from 1 to 15, and we have $\delta(T)=(3,2,0,2,0,2,0,0,0,3,0,0,2,0,0)\in \mathcal{B}_{14,15}$.
\end{exam}

Lemma \ref{planeB} has a fairly straightforward proof by induction and will be omitted here. There is a more general result that plane forests are in bijection with a sequence of unit compositions. More details can be found in \cite[p.34]{EC2}, in which the author use the terminology of {\em {\L}ukasiewicz words} instead of compositions. In the remainder of this section we concentrate on plane trees with a specified vertex of given outdegree. Now we are ready to establish our bijection.

\begin{theo}\label{bij_plane}
There is a bijection between $\mathcal{T}_{n,i}$ and $\mathcal{A}_{n-i,n}$.
\end{theo}
\pf Given an ordered pair $(T,v) \in \mathcal{T}_{n,i}$, we label the vertices of $T$ in depth-first order. Suppose $v$ is labeled $j$, $1 \leq j \leq n+1$. We define
\[\bar{\delta}(T,v)=\alpha=(d_{j+1},\ldots, d_{n+1}, d_1, \ldots, d_{j-1}).\]
For example, let $T$ be the plane tree shown in Figure \ref{Fig_plane}, and the circled vertex with label $4$ is the specified vertex $v$, which has outdegree $2$. Then we have $\bar{\delta}(T,v)=(0,2,0,0,0,3,0,0,2,0,0,3,2,0).$
Since there are $n$ edges of $T$, and $v$ has outdegree $i$, it is obvious that $\sum_{t=1,t\neq j}^{n+1}d_{t}=n-i$, thus $\bar{\delta}(T,v) \in \mathcal{A}_{n-i,n}$.

On the other hand, given $\alpha \in \mathcal{A}_{n-i,n}$, we can uniquely decompose $\alpha$ into the form$$\alpha=\alpha_1 \alpha_2 \cdots \alpha_s \alpha_0$$
for some nonnegative integer $s$ such that $\alpha_1,\alpha_2, \ldots, \alpha_s$ are unit compositions, and $\alpha_0$ is a positive composition with $f(\alpha_0)=s-i$. We call such a decomposition the {\em fundamental decomposition} of $\alpha$. For example, the fundamental decomposition of $\alpha=(0,2,0,0,0,3,0,0,2,0,0,3,2,0)$ is the following:
\[(0)~(2,0,0)~(0)~(3,0,0,2,0,0)~(3,2,0).\]
Now we set $\alpha^{\prime}=\alpha_0 ~i~ \alpha_1 \alpha_2 \cdots \alpha_s$. One can easily verify that $\alpha^{\prime}$ is a unit composition. We set $T=\delta^{-1}(\alpha^{\prime})$. Suppose $\alpha_0$ is of length $l$. Let $v$ be the $(l+1)$-th vertex of $T$ in the depth-first order.
Then we have $(T,v)=\bar{\delta}^{-1}(\alpha)$. Hence we proved that the map $\bar{\delta}$ is a bijection.
\qed

\noindent{\em Proof of Theorem \ref{main_plane}}:
It is a basic result in enumerative combinatorics with an easy bijective proof that the number of $k$ compositions of $n$, or equivalently, the number of nonnegative integer solutions of the equation $x_1+x_2+\cdots+x_k=n$ is ${n+k-1 \choose n}$, so our result then follows.
\qed

Given a sequence of nonnegative integers $(r_0,r_1, \ldots, r_n)$ with $\sum_{j=0}^{n}r_j=n+1$ and $\sum_{j=0}^{n}jr_j=n$. The number of plane trees with $n+1$ vertices such that exactly $r_j$ vertices have outdegree $j$ is given by
\begin{equation}\label{ECplanetree}
\frac{1}{n+1}{n+1 \choose r_0, r_1, \dots, r_n}.
\end{equation}
 A more general version of the above result concerning the enumeration of plane forests (graphs
such that every connected component is a plane tree) by outdegree sequence is given in \cite[Theorem 5.3.10]{EC2}. 

Combining Theorem \ref{main_plane} and formula \eqref{ECplanetree} we have the following identity:
 \begin{equation}
 \sum\frac{r_i}{n+1}{n+1 \choose r_0, r_1, \ldots, r_n}={2n-i-1 \choose n-1},
 \end{equation}
 where the sum is over  all sequences of nonnegative integers $(r_0,r_1, \ldots, r_n)$ such that  $\sum_{j=0}^{n}r_j=n+1$ and $\sum_{j=0}^{n}jr_j=n$.

\section{Counting vertices with given degree in plane trees}

In this section we count vertices among all $n$-plane trees with given degree. Note that for each vertex $v$ in a plane tree with outdegree $i$, the degree of $v$ is also $i$ when $v$ is the root, otherwise the degree of $v$ is $i+1$. Thus it is sufficient to count the number of $n$-plane trees with specified root degree, and the number of vertices with given outdegree among all $n$-plane trees with specified root degree. 

\begin{lem}\label{lem_root_m}
For positive integers $n$ and $m$ with $m \leq n$, the total number of $n$-plane trees with root degree $m$ is 
\begin{equation}\label{plane_root_m}
\frac{m}{n}{2n-m-1 \choose n-1}.
\end{equation}
\end{lem}
\pf Let $T$ be an $n$-plane tree with root degree $m$. We have $\delta(T)=(m, d_2, d_3, \ldots, d_{n+1})$. The total number of nonegative integer solutions of the equations $d_2+d_3+\cdots+d_{n+1}=n-m$ is ${2n-m-1 \choose n-1}$. On the other hand, from the discussion in the previous section we know that among the $n$ cyclic permutations of $d_2, d_3, \ldots, d_{n+1}$, exactly $m$ of them are valid outdegree sequences, therefore \eqref{plane_root_m} holds.
\qed

\begin{theo}\label{rootdeg}
For integers $n,m \geq 1$ and $i \geq 0$, the total number of vertices of outdegree $i$ among all $n$-plane trees with root degree $m$ is 
\begin{numcases}{}
m {2n-m-i-2 \choose n-2}, \mbox{when\ }  i \neq m, \label{in=m}\\
m {2n-2m-2 \choose n-2}+\frac{m}{n}{2n-m-1 \choose n-1},  \mbox{when\ } i = m.\label{i=m}
\end{numcases}
\end{theo}

\pf We first consider the case when $i \neq m$. Let $T$ be an $n$-plane tree with root degree $m$, and $v$ be a specified vertex in $T$ with outdegree $i$. We have that $m$ appears at least once in $\bar{\delta}(T,v)$. The total number of nonegative integer solutions of the equations 
\begin{equation}\label{eq_n-m-i}
x_1+x_2+\cdots+x_{n-1}=n-m-i
\end{equation}
is ${2n-m-i-2 \choose n-2}$. 
For any nonnegative integer sequence $x_1, x_2, \ldots, x_{n-1}$ satisfying equation \eqref{eq_n-m-i},  it can be uniquely decomposed into $a_1 a_2 \cdots a_s b_0$ such that $a_1, a_2, \ldots, a_s$ are unit compositions and $b_0$ is a positive composition such that $f(b_0)=s-m-i+1$. There are $m$ ways to insert  $m$ and $i$ to form a sequence $\alpha$ such that $\bar{\delta}^{-1}(\alpha)=(T,v)$, where $T$ is an $n$-plane tree with root degree $m$ and $v$ s a specified vertex of $T$ with outdegree $i$:
\begin{eqnarray*}
& b_0 \ i \ a_1\cdots a_{s-m+1}\ a_{s-m+2}\cdots a_{s-1} \ a_s \\
& a_s \ b_0 \ i \  a_1\cdots a_{s-m+1}\ a_{s-m+2} \cdots a_{s-1}\\
& \cdots \\
& a_{s-m+2} \cdots a_{s-1}\  a_{s}\ b_0 \ i\  a_1\cdots a_{s-m+1}
\end{eqnarray*}
Therefore \eqref{in=m} holds.
For the case $i=m$, the roots should be also counted. Combining \eqref{in=m} and \eqref{plane_root_m} we get \eqref{i=m}.
\qed

Combining the above result with Theorem \ref{main_plane} we get Corollary \ref{coro_plane}.

\noindent {\em Proof of Corollary \ref{coro_plane}:}
Among all $n$-plane trees, the total number of vertices of degree $i$ equals the total number of vertices of outdegree $i-1$ minus the number of $n$-plane trees with root degree $i-1$ plus the number of $n$-plane trees with root degree $i$. Corollary \ref{coro_plane} then follows by simple calculation. \qed

Note that more recently Eu, Seo and Shin \cite{EuSeoShin} enumerated vertices among all plane trees with levels and degrees, and obtained various interesting results. 
Theorem \ref{main_plane} and Corollary \ref{coro_plane} can also be obtained from their results by a telescoping summation. 

\section{Counting vertices with given outdegree in $k$-ary trees}


Let $\mathcal{T}^k_{n}$ denote the set of $k$-ary trees with $n$ edges, and $\mathcal{T}^k_{n,i}$ denote the set of ordered pairs $(T,v)$ such that $T\in \mathcal{T}^k_n$, and $v$ is a vertex of $T$ of outdegree $i$. Given $(T,v)\in \mathcal{T}^k_{n,i}$, let $T^{\prime}$ denote the complete $k$-ary tree that corresponds to $T$. Note that $v$ becomes an internal vertex in $T^{\prime}$ (and has outdegree $k$), which we denote as $v^{\prime}$. Let $\alpha=\bar{\delta}(T^{\prime},v^{\prime})$. Since there are $n+1$ internal vertices and $k(n+1)$ edges in $T'$, and each internal vertices of $T^{\prime}$ has outdegree $k$, from Theorem \ref{bij_plane} we have the following corollary.

\begin{coro}\label{kary2comp}
The map $\bar{\delta}$ is a bijection between $\mathcal{T}^k_{n,i}$ and compositions $\alpha$ that satisfy the following:
\begin{itemize}
\item[(i)] $\alpha \in \mathcal{A}_{kn,k(n+1)}$, and $\alpha$ consists of exactly $n$ $k$'s and $(kn+k-n)$ $0$'s;
\item[(ii)] In the fundamental decomposition $\alpha=\alpha_1 \alpha_2 \cdots \alpha_s \alpha_0$, we have $s\geq k$, and among the first $k$ unit compositions $\alpha_1, \alpha_2, \ldots, \alpha_k$,  exactly $i$ of them begin with $k$.
\end{itemize}
\end{coro}

\begin{exam}\label{examp_kary}
Let $k=3$, $n=8$ and $i=2$. Figure \ref{Figkary} shows a ternary tree $T$ with $8$ edges and the corresponding complete ternary tree $T^{\prime}$ with $k(n+1)+1=28$ vertices. A specified vertex $v$ of outdegree 2 is circled in $T$, and which corresponds to an internal vertex $v'$ in $T'$. We have $$\alpha=\bar{\delta}(T^{\prime},v^{\prime})=(3,0,0,0,0,3,0,0,0,0,3,0,0,3,3,0,0,0,3,0,0,0,0,3,3,0,0).$$
The fundamental decomposition of $\alpha$ is
$$(3,0,0,0)~(0)~(3,0,0,0)~(0)~(3,0,0,3,3,0,0,0,3,0,0,0,0)~(3,3,0,0).$$
It is easy to check that $\alpha$ satisfies conditions (i) and (ii).

\begin{figure}
\begin{center}
\begin{picture}(400,100)
\put(0,0){
\begin{picture}(130,100)
\setlength{\unitlength}{12pt}
\thicklines
\put(5,6){\line(-3,-2){3}}
\put(5,6){\line(3,-2){3}}
\put(2,4){\line(1,-1){2}}
\put(4,2){\line(-1,-1){2}}
\put(4,2){\line(1,-1){2}}
\put(8,4){\line(1,-1){2}}
\put(10,2){\line(-1,-1){2}}
\put(10,2){\line(0,-1){2}}
\thinlines
\put(5,6){\circle*{0.25}}
\put(2,4){\circle*{0.25}}
\put(4,2){\circle*{0.25}}\put(4,2){\color{red}\circle{0.5}}\put(4.3,2){\small $v$}
\put(2,0){\circle*{0.25}}
\put(6,0){\circle*{0.25}}
\put(8,4){\circle*{0.25}}
\put(10,2){\circle*{0.25}}
\put(8,0){\circle*{0.25}}
\put(10,0){\circle*{0.25}}

\put(3,6){$T$}
\end{picture}}

\put(200,0){
\begin{picture}(200,100)
\setlength{\unitlength}{16pt}
\thicklines
\put(5,6){\line(-2,-1){4}}
\put(5,6){\line(2,-1){4}}
\put(1,4){\line(1,-1){2}}
\put(3,2){\line(-1,-1){2}}
\put(3,2){\line(1,-1){2}}
\put(9,4){\line(1,-1){2}}
\put(11,2){\line(-3,-2){3}}
\put(11,2){\line(0,-1){2}}
\thinlines
\put(5,6){\circle*{0.2}}\put(4.7,6.2){\small $1$}
\put(1,4){\circle*{0.2}}\put(0.6,4){\small $2$}
\put(3,2){\circle*{0.2}}\put(3,2){{\bf \color{red}\circle{0.4}}}\put(3.3,2){\small $5$($v^{\prime}$)}
\put(1,0){\circle*{0.2}}\put(0.4,0){\small $6$}
\put(5,0){\circle*{0.2}}\put(5.2,0){\small $11$}
\put(9,4){\circle*{0.2}}\put(9.2,4){\small $16$}
\put(11,2){\circle*{0.2}}\put(11.2,2){\small $19$}
\put(8,0){\circle*{0.2}}\put(7.2,0){\small $20$}
\put(11,0){\circle*{0.2}}\put(10.2,0){\small $24$}
\put(5,6){{\bf \color{blue}\line(0,-1){2}}}
\put(1,4){{\bf \color{blue}\line(-1,-1){2}}}
\put(1,4){{\bf \color{blue}\line(0,-1){2}}}
\put(3,2){{\bf \color{blue}\line(0,-1){2}}}
\put(9,4){{\bf \color{blue}\line(-1,-1){2}}}
\put(9,4){{\bf \color{blue}\line(0,-1){2}}}
\put(11,2){{\bf \color{blue}\line(1,-2){1}}}
\put(1,0){{\bf \color{blue}\line(-1,-2){1}}}
\put(1,0){{\bf \color{blue}\line(0,-1){2}}}
\put(1,0){{\bf \color{blue}\line(1,-2){1}}}
\put(5,0){{\bf \color{blue}\line(-1,-2){1}}}
\put(5,0){{\bf \color{blue}\line(0,-1){2}}}
\put(5,0){{\bf \color{blue}\line(1,-2){1}}}
\put(8,0){{\bf \color{blue}\line(-1,-2){1}}}
\put(8,0){{\bf \color{blue}\line(0,-1){2}}}
\put(8,0){{\bf \color{blue}\line(1,-2){1}}}
\put(11,0){{\bf \color{blue}\line(-1,-2){1}}}
\put(11,0){{\bf \color{blue}\line(0,-1){2}}}
\put(11,0){{\bf \color{blue}\line(1,-2){1}}}

\put(-1,2){{\bf \color{blue}\circle*{0.1}}}\put(-1.4,2){\small $3$}
\put(1,2){{\bf \color{blue}\circle*{0.1}}}\put(1.2,2){\small $4$}
\put(0,-2){{\bf \color{blue}\circle*{0.1}}}\put(-0.2,-2.5){\small $7$}
\put(1,-2){\bf \color{blue}{\circle*{0.1}}}\put(0.8,-2.5){\small $8$}
\put(2,-2){\bf \color{blue}{\circle*{0.1}}}\put(1.8,-2.5){\small $9$}
\put(3,0){\bf \color{blue}{\circle*{0.1}}}\put(3.2,0){\small $10$}
\put(4,-2){\bf \color{blue}{\circle*{0.1}}}\put(3.8,-2.5){\small $12$}
\put(5,-2){\bf \color{blue}{\circle*{0.1}}}\put(4.8,-2.5){\small $13$}
\put(6,-2){\bf \color{blue}{\circle*{0.1}}}\put(5.8,-2.5){\small $14$}
\put(5,4){\bf \color{blue}{\circle*{0.1}}}\put(5.2,4){\small $15$}
\put(7,2){\bf \color{blue}{\circle*{0.1}}}\put(6.2,2){\small $17$}
\put(9,2){\bf \color{blue}{\circle*{0.1}}}\put(9.2,2){\small $18$}
\put(7,-2){\bf \color{blue}{\circle*{0.1}}}\put(6.8,-2.5){\small $21$}
\put(8,-2){\bf \color{blue}{\circle*{0.1}}}\put(7.8,-2.5){\small $22$}
\put(9,-2){\bf \color{blue}{\circle*{0.1}}}\put(8.8,-2.5){\small $23$}
\put(10,-2){\bf \color{blue}{\circle*{0.1}}}\put(9.8,-2.5){\small $25$}
\put(11,-2){\bf \color{blue}{\circle*{0.1}}}\put(10.8,-2.5){\small $26$}
\put(12,-2){\bf \color{blue}{\circle*{0.1}}}\put(11.8,-2.5){\small $27$}
\put(12,0){\bf \color{blue}{\circle*{0.1}}}\put(12.2,0){\small $28$}

\put(3,6){$T^{\prime}$}
\end{picture}}
\end{picture}
\vskip 20mm
\caption{A ternary tree $T$ and the corresponding complete ternary tree $T^{\prime}$.}
\label{Figkary}
\end{center}
\end{figure}
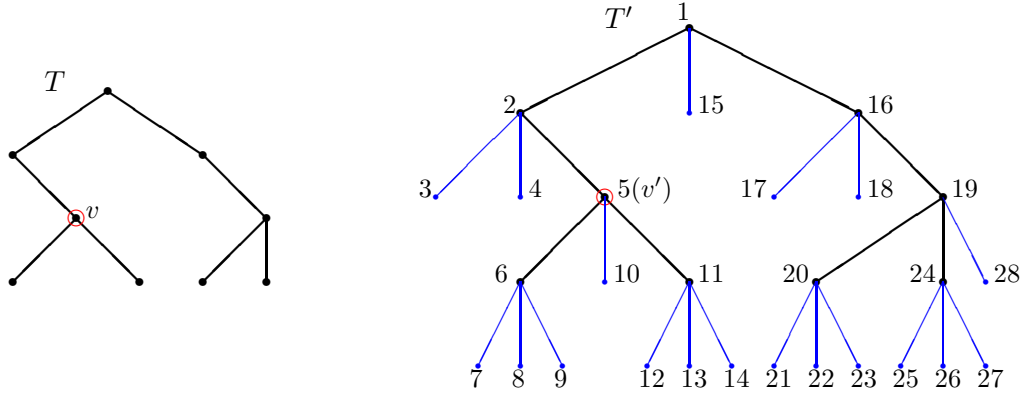

\end{exam}

Let $\mathcal{A}^*_{kn,k(n+1)}$ denote the set of compositions in $\mathcal{A}_{kn,k(n+1)}$ that satisfy (i) and (ii).
From Corollary \ref{kary2comp} we know that to prove Theorem \ref{main_kary}, we need to show the following.
\begin{theo}\label{bijXY}
For integers $n,k \geq 1$ and $i \geq 0$,
there is a bijection between $\mathcal{A}^*_{kn,k(n+1)}$ and the set of
of ordered pairs $(X,Y)$, where $X$ is an $i$-element subset of $[k]:=\{1,2, \ldots, k\}$, and $Y$ is an $(n-i)$-element subset of $[kn]$.
\end{theo}
\pf Given $\alpha \in \mathcal{A}^*_{kn,k(n+1)}$. Suppose the fundamental decomposition of $\alpha$ is $\alpha=\alpha_1 \alpha_2 \cdots \alpha_s \alpha_0$. Let $\alpha_{j_1}, \alpha_{j_2}, \ldots, \alpha_{j_i}$ be the $i$ unit compositions among $\alpha_1, \alpha_2, \ldots, \alpha_k$ that begin with $k$. Set $X=\{j_1, j_2, \ldots, j_i\}$. Let $\beta=(b_1, b_2, \ldots, b_{kn})$ be the composition obtained from $\alpha$ by deleting the first number in each of the first $k$ unit compositions $\alpha_i$, $1 \leq i \leq k$. Then exactly $n-i$ members among $b_1, b_2, \ldots, b_{kn}$ equals $k$. Suppose they are $b_{l_1}, b_{l_2}, \ldots, b_{l_{n-i}}$, with $l_1, l_2, \cdots, l_{n-i} \in [kn]$. Now we set $Y=\{l_1, l_2, \ldots, l_{n-i}\}$ and $\phi(\alpha)=(X,Y)$.

For example, for the composition $\alpha$ in Example \ref{examp_kary}, we have $\phi(\alpha)=(X,Y)$ where $X=\{1,3\}$, and $Y=\{8,11,12,16,21,22\}$.

Now we prove that $\phi$ is a bijection by defining its inverse. Given $(X,Y)$ with $X=\{j_1, j_2, \ldots, j_i\} \subseteq [k]$ and $Y=\{l_1, l_2, \ldots, l_{n-i}\} \subseteq [kn]$. we define a sequence $\beta=(b_1,b_2, \ldots, b_{kn})$ such that $b_j=k$ if $j \in Y$, and $b_j=0$ otherwise. Similarly we define $\gamma=(c_1,c_2, \ldots, c_{k})$ to be a sequence of integers such that $c_j=k$ if $j \in X$, and $c_j=0$ otherwise. Now we insert $c_1,c_2, \ldots, c_{k}$ into $\beta$ one-by-one to get a sequence $\alpha=(a_1, a_2, \ldots, a_{k(n+1)})$ such that in the fundamental decomposition of $\alpha$, the first $k$ unit compositions begin with $c_1,c_2, \ldots, c_{k}$, respectively. Note that when inserting a $0$, it forms a unit composition by itself; when inserting a $k$, it will ``use" $k$ $0$'s from $\beta$, since there are $i$ $k$'s among $c_1,c_2, \ldots, c_{k}$ and $kn-k(n-i)=ki$ ``remaining" $0$'s (each of the $n-i$ $k$'s in $\beta$ will ``use" $k$ $0$'s too), such an insertion is always possible. Thus we proved that $\phi$ is a bijection.\qed

Theorem \ref{main_kary} follows immediately from Theorem \ref{bijXY}.

\begin{exam}
Let $n=2$, $k=2$ and $i=1$. There are eight pairs $(X,Y)$, where $X$ is a one-element subset of $\{1,2\}$, and $Y$ is a one-element subset of $\{1,2,3,4\}$. Table \ref{tabbinary8}
shows the corresponding $\alpha$, $(T',v')$ and $(T,v)$ for each pair $(X,Y)$.
\end{exam}

\begin{table}\centering
\begin{tabular}{|c|c|c|c|c||c|c|c|c|c|}\hline
$X$&$Y$&$\alpha$&($T'$,$v'$)&($T$,$v$)&$X$&$Y$&$\alpha$&($T'$,$v'$)&($T$,$v$)\\ \hline
\raisebox{15pt}[0pt]{\{1\}}&\raisebox{15pt}[0pt]{\{1\}}&\raisebox{15pt}[0pt]{({\bf \color{red}2}2000)({\bf \color{red}0})}&
\begin{picture}(40,40)
\setlength{\unitlength}{10pt}
\thicklines
\put(3,3){\line(-1,-1){1}}
\put(2,2){\line(-1,-1){1}}
\thinlines
\put(3,3){\circle*{0.2}}\put(3,3){{\bf \color{red}\circle{0.4}}}
\put(2,2){\circle*{0.2}}
\put(1,1){\circle*{0.2}}
\put(3,3){{\bf \color{blue}\line(1,-1){1}}}\put(4,2){{\bf \color{blue}\circle*{0.1}}}
\put(2,2){{\bf \color{blue}\line(1,-1){1}}}\put(3,1){{\bf \color{blue}\circle*{0.1}}}
\put(1,1){{\bf \color{blue}\line(-1,-1){1}}}\put(0,0){{\bf \color{blue}\circle*{0.1}}}
\put(1,1){{\bf \color{blue}\line(1,-1){1}}}\put(2,0){{\bf \color{blue}\circle*{0.1}}}
\end{picture}&
\begin{picture}(20,40)
\setlength{\unitlength}{10pt}
\thicklines
\put(2,3){\line(-1,-1){1}}
\put(1,2){\line(-1,-1){1}}
\thinlines
\put(2,3){\circle*{0.2}}\put(2,3){{\bf \color{red}\circle{0.4}}}
\put(1,2){\circle*{0.2}}
\put(0,1){\circle*{0.2}}
\end{picture}&
\raisebox{15pt}[0pt]{\{2\}}&\raisebox{15pt}[0pt]{\{1\}}&\raisebox{15pt}[0pt]{({\bf \color{red}0})({\bf \color{red}2}2000)}&
\begin{picture}(40,40)
\setlength{\unitlength}{10pt}
\thicklines
\put(2,3){\line(1,-1){1}}
\put(3,2){\line(-1,-1){1}}
\thinlines
\put(2,3){\circle*{0.2}}\put(2,3){{\bf \color{red}\circle{0.4}}}
\put(3,2){\circle*{0.2}}
\put(2,1){\circle*{0.2}}
\put(2,3){{\bf \color{blue}\line(-1,-1){1}}}\put(1,2){{\bf \color{blue}\circle*{0.1}}}
\put(3,2){{\bf \color{blue}\line(1,-1){1}}}\put(4,1){{\bf \color{blue}\circle*{0.1}}}
\put(2,1){{\bf \color{blue}\line(-1,-1){1}}}\put(1,0){{\bf \color{blue}\circle*{0.1}}}
\put(2,1){{\bf \color{blue}\line(1,-1){1}}}\put(3,0){{\bf \color{blue}\circle*{0.1}}}
\end{picture}&
\begin{picture}(20,40)
\setlength{\unitlength}{10pt}
\thicklines
\put(1,3){\line(1,-1){1}}
\put(2,2){\line(-1,-1){1}}
\thinlines
\put(1,3){\circle*{0.2}}\put(1,3){{\bf \color{red}\circle{0.4}}}
\put(2,2){\circle*{0.2}}
\put(1,1){\circle*{0.2}}
\end{picture}
\\ \hline
\raisebox{15pt}[0pt]{\{1\}}&\raisebox{15pt}[0pt]{\{2\}}&\raisebox{15pt}[0pt]{({\bf \color{red}2}0200)({\bf \color{red}0})}&
\begin{picture}(40,40)
\setlength{\unitlength}{10pt}
\thicklines
\put(2.5,3){\line(-1,-1){1}}
\put(1.5,2){\line(1,-1){1}}
\thinlines
\put(2.5,3){\circle*{0.2}}\put(2.5,3){{\bf \color{red}\circle{0.4}}}
\put(1.5,2){\circle*{0.2}}
\put(2.5,1){\circle*{0.2}}
\put(2.5,3){{\bf \color{blue}\line(1,-1){1}}}\put(3.5,2){{\bf \color{blue}\circle*{0.1}}}
\put(1.5,2){{\bf \color{blue}\line(-1,-1){1}}}\put(0.5,1){{\bf \color{blue}\circle*{0.1}}}
\put(2.5,1){{\bf \color{blue}\line(-1,-1){1}}}\put(1.5,0){{\bf \color{blue}\circle*{0.1}}}
\put(2.5,1){{\bf \color{blue}\line(1,-1){1}}}\put(3.5,0){{\bf \color{blue}\circle*{0.1}}}
\end{picture}
&
\begin{picture}(20,40)
\setlength{\unitlength}{10pt}
\thicklines
\put(1.5,3){\line(-1,-1){1}}
\put(0.5,2){\line(1,-1){1}}
\thinlines
\put(1.5,3){\circle*{0.2}}\put(1.5,3){{\bf \color{red}\circle{0.4}}}
\put(0.5,2){\circle*{0.2}}
\put(1.5,1){\circle*{0.2}}
\end{picture}
&\raisebox{15pt}[0pt]{\{2\}}&\raisebox{15pt}[0pt]{\{2\}}&\raisebox{15pt}[0pt]{({\bf \color{red}0})({\bf \color{red}2}0200)}&
\begin{picture}(40,40)
\setlength{\unitlength}{10pt}
\thicklines
\put(1,3){\line(1,-1){1}}
\put(2,2){\line(1,-1){1}}
\thinlines
\put(1,3){\circle*{0.2}}\put(1,3){{\bf \color{red}\circle{0.4}}}
\put(2,2){\circle*{0.2}}
\put(3,1){\circle*{0.2}}
\put(1,3){{\bf \color{blue}\line(-1,-1){1}}}\put(0,2){{\bf \color{blue}\circle*{0.1}}}
\put(2,2){{\bf \color{blue}\line(-1,-1){1}}}\put(1,1){{\bf \color{blue}\circle*{0.1}}}
\put(3,1){{\bf \color{blue}\line(-1,-1){1}}}\put(2,0){{\bf \color{blue}\circle*{0.1}}}
\put(3,1){{\bf \color{blue}\line(1,-1){1}}}\put(4,0){{\bf \color{blue}\circle*{0.1}}}
\end{picture}
&\begin{picture}(20,40)
\setlength{\unitlength}{10pt}
\thicklines
\put(0,3){\line(1,-1){1}}
\put(1,2){\line(1,-1){1}}
\thinlines
\put(0,3){\circle*{0.2}}\put(0,3){{\bf \color{red}\circle{0.4}}}
\put(1,2){\circle*{0.2}}
\put(2,1){\circle*{0.2}}
\end{picture}
\\ \hline
\raisebox{15pt}[0pt]{\{1\}}&\raisebox{15pt}[0pt]{\{3\}}&\raisebox{15pt}[0pt]{({\bf \color{red}2}00)({\bf \color{red}0})(20)}&
\begin{picture}(40,40)
\setlength{\unitlength}{10pt}
\thicklines
\put(2,3){\line(1,-1){1}}
\put(3,2){\line(-1,-1){1}}
\thinlines
\put(2,3){\circle*{0.2}}
\put(3,2){\circle*{0.2}}\put(3,2){{\bf \color{red}\circle{0.4}}}
\put(2,1){\circle*{0.2}}
\put(2,3){{\bf \color{blue}\line(-1,-1){1}}}\put(1,2){{\bf \color{blue}\circle*{0.1}}}
\put(3,2){{\bf \color{blue}\line(1,-1){1}}}\put(4,1){{\bf \color{blue}\circle*{0.1}}}
\put(2,1){{\bf \color{blue}\line(-1,-1){1}}}\put(1,0){{\bf \color{blue}\circle*{0.1}}}
\put(2,1){{\bf \color{blue}\line(1,-1){1}}}\put(3,0){{\bf \color{blue}\circle*{0.1}}}
\end{picture}
&\begin{picture}(20,40)
\setlength{\unitlength}{10pt}
\thicklines
\put(1,3){\line(1,-1){1}}
\put(2,2){\line(-1,-1){1}}
\thinlines
\put(1,3){\circle*{0.2}}
\put(2,2){\circle*{0.2}}\put(2,2){{\bf \color{red}\circle{0.4}}}
\put(1,1){\circle*{0.2}}
\end{picture}
&\raisebox{15pt}[0pt]{\{2\}}&\raisebox{15pt}[0pt]{\{3\}}&\raisebox{15pt}[0pt]{({\bf \color{red}0})({\bf \color{red}2}0020)}&
\begin{picture}(40,40)
\setlength{\unitlength}{10pt}
\thicklines
\put(1,3){\line(1,-1){1}}
\put(2,2){\line(1,-1){1}}
\thinlines
\put(1,3){\circle*{0.2}}
\put(2,2){\circle*{0.2}}\put(2,2){{\bf \color{red}\circle{0.4}}}
\put(3,1){\circle*{0.2}}
\put(1,3){{\bf \color{blue}\line(-1,-1){1}}}\put(0,2){{\bf \color{blue}\circle*{0.1}}}
\put(2,2){{\bf \color{blue}\line(-1,-1){1}}}\put(1,1){{\bf \color{blue}\circle*{0.1}}}
\put(3,1){{\bf \color{blue}\line(-1,-1){1}}}\put(2,0){{\bf \color{blue}\circle*{0.1}}}
\put(3,1){{\bf \color{blue}\line(1,-1){1}}}\put(4,0){{\bf \color{blue}\circle*{0.1}}}
\end{picture}
&\begin{picture}(20,40)
\setlength{\unitlength}{10pt}
\thicklines
\put(0,3){\line(1,-1){1}}
\put(1,2){\line(1,-1){1}}
\thinlines
\put(0,3){\circle*{0.2}}
\put(1,2){\circle*{0.2}}\put(1,2){{\bf \color{red}\circle{0.4}}}
\put(2,1){\circle*{0.2}}
\end{picture}
\\ \hline
\raisebox{15pt}[0pt]{\{1\}}&\raisebox{15pt}[0pt]{\{4\}}&\raisebox{15pt}[0pt]{({\bf \color{red}2}00)({\bf \color{red}0})(0)(2)}&
\begin{picture}(40,40)
\setlength{\unitlength}{10pt}
\thicklines
\put(3,3){\line(-1,-1){1}}
\put(2,2){\line(-1,-1){1}}
\thinlines
\put(3,3){\circle*{0.2}}
\put(2,2){\circle*{0.2}}\put(2,2){{\bf \color{red}\circle{0.4}}}
\put(1,1){\circle*{0.2}}
\put(3,3){{\bf \color{blue}\line(1,-1){1}}}\put(4,2){{\bf \color{blue}\circle*{0.1}}}
\put(2,2){{\bf \color{blue}\line(1,-1){1}}}\put(3,1){{\bf \color{blue}\circle*{0.1}}}
\put(1,1){{\bf \color{blue}\line(-1,-1){1}}}\put(0,0){{\bf \color{blue}\circle*{0.1}}}
\put(1,1){{\bf \color{blue}\line(1,-1){1}}}\put(2,0){{\bf \color{blue}\circle*{0.1}}}
\end{picture}
&\begin{picture}(20,40)
\setlength{\unitlength}{10pt}
\thicklines
\put(2,3){\line(-1,-1){1}}
\put(1,2){\line(-1,-1){1}}
\thinlines
\put(2,3){\circle*{0.2}}
\put(1,2){\circle*{0.2}}\put(1,2){{\bf \color{red}\circle{0.4}}}
\put(0,1){\circle*{0.2}}
\end{picture}
&\raisebox{15pt}[0pt]{\{2\}}&\raisebox{15pt}[0pt]{\{4\}}&\raisebox{15pt}[0pt]{({\bf \color{red}0})({\bf \color{red}2})00(0)(2)}&
\begin{picture}(40,40)
\setlength{\unitlength}{10pt}
\thicklines
\put(2.5,3){\line(-1,-1){1}}
\put(1.5,2){\line(1,-1){1}}
\thinlines
\put(2.5,3){\circle*{0.2}}
\put(1.5,2){\circle*{0.2}}\put(1.5,2){{\bf \color{red}\circle{0.4}}}
\put(2.5,1){\circle*{0.2}}
\put(2.5,3){{\bf \color{blue}\line(1,-1){1}}}\put(3.5,2){{\bf \color{blue}\circle*{0.1}}}
\put(1.5,2){{\bf \color{blue}\line(-1,-1){1}}}\put(0.5,1){{\bf \color{blue}\circle*{0.1}}}
\put(2.5,1){{\bf \color{blue}\line(-1,-1){1}}}\put(1.5,0){{\bf \color{blue}\circle*{0.1}}}
\put(2.5,1){{\bf \color{blue}\line(1,-1){1}}}\put(3.5,0){{\bf \color{blue}\circle*{0.1}}}
\end{picture}
&\begin{picture}(20,40)
\setlength{\unitlength}{10pt}
\thicklines
\put(1.5,3){\line(-1,-1){1}}
\put(0.5,2){\line(1,-1){1}}
\thinlines
\put(1.5,3){\circle*{0.2}}
\put(0.5,2){\circle*{0.2}}\put(0.5,2){{\bf \color{red}\circle{0.4}}}
\put(1.5,1){\circle*{0.2}}
\end{picture}
\\ \hline
\end{tabular}
\caption{The bijection for 8 pairs $(T,v)$ when $n=2$, $k=2$ and $i=1$.}\label{tabbinary8}
\end{table}

\vskip 3cm

\vskip 2mm \noindent{\bf Acknowledgments.} 
This work is partially supported by the Science and Technology Commission of Shanghai Municipality (STCSM), No. 13dz2260400.

\end{document}